\documentclass{amsart}

\usepackage{hyperref}
\usepackage{enumerate}
\usepackage{rotating}
\usepackage[matrix,arrow,graph,curve,frame]{xy}

\theoremstyle{plain}      
    \newtheorem{theorem}{Theorem}[section]
    \newtheorem{proposition}{Proposition}[section]

\theoremstyle{definition}
    
    \newtheorem{example}{Example}[section]
\theoremstyle{remark}

\usepackage{mathrsfs}\renewcommand{\mathcal}{\mathscr}
\newcommand{\A}{\ensuremath{\mathcal{A}}}
\newcommand{\B}{\ensuremath{\mathcal{B}}}
\newcommand{\C}{\ensuremath{\mathcal{C}}}
\newcommand{\D}{\ensuremath{\mathcal{D}}}

\newcommand{\copr}{\ensuremath{\mathrm{copr}}}

\newcommand{\Comod}{\ensuremath{\mathbf{Comod}}}
\newcommand{\Co}{\ensuremath{\mathbf{Co}}}
\newcommand{\Lex}{\ensuremath{\mathrm{Lex}}}

\newcommand{\Set}{\ensuremath{\mathbf{Set}}}
\newcommand{\Vect}{\ensuremath{\mathbf{Vect}}}
\newcommand{\Endv}{\ensuremath{\mathrm{End}^\vee}}
\newcommand{\dint}{\ensuremath{{\displaystyle \int}}}
\newcommand{\tint}{\ensuremath{{\textstyle \int}}}

\newcommand{\ox}{\ensuremath{\otimes}}
\newcommand{\op}{\ensuremath{\mathrm{op}}}

\newcommand{\ob}{\ensuremath{\mathrm{ob}}}

\newcommand{\ra}{\ensuremath{\xymatrix@=20pt@1{\ar[r]&}}}

\begin{document}
\title{On endomorphism algebras of separable monoidal functors}
\author{Brian Day}
\author{Craig Pastro}
\thanks{The first author gratefully acknowledges partial support of an
Australian Research Council grant while the second gratefully acknowledges
support of an international Macquarie University Research Scholarship.}
\address{Department of Mathematics \\
         Macquarie University \\
         New South Wales 2109 Australia}
\email{craig@ics.mq.edu.au}
\date{\today}

\begin{abstract}
We show that the (co)endomorphism algebra of a sufficiently separable
``fibre'' functor into $\Vect_k$, for $k$ a field of characteristic $0$,
has the structure of what we call a ``unital'' von Neumann core in
$\Vect_k$. For $\Vect_k$, this particular notion of algebra is weaker than
that of a Hopf algebra, although the corresponding concept in $\Set$ is
again that of a group.
\end{abstract}
\maketitle

\section{Introduction}

Let $(\C,\ox,I,c)$ be a symmetric (or just braided) monoidal category. Recall
that an \emph{algebra} in $\C$ is an object $A \in \C$ equipped with a
multiplication $\mu:A \ox A \ra A$ and a unit $\eta:I \ra A$ satisfying
$\mu_3 = \mu(1 \ox \mu) = \mu(\mu \ox 1):A^{\ox 3} \ra A$ (associativity)
and $\mu(\eta \ox 1) = 1 = \mu(1 \ox \eta):A \ra A$ (unit conditions).
Dually, a \emph{coalgebra} in $\C$ is an object $C \in \C$ equipped with a
comultiplication $\delta:C \ra C \ox C$ and a counit $\epsilon:C \ra I$
satisfying $\delta_3 = (1 \ox \delta)\delta = (\delta \ox 1)\delta: C \ra
C^{\ox 3}$ (coassociativity) and $(\epsilon \ox 1)\delta = 1 = (1 \ox
\epsilon)\delta:C \ra C$ (counit conditions).

A \emph{very weak bialgebra} in $\C$ is an object $A \in \C$ with both the
structure of an algebra and a coalgebra in $\C$ related by the axiom
\[
  \delta\mu=(\mu \ox \mu)(1 \ox c \ox 1)(\delta \ox \delta):A \ox A \ra A \ox A.
\]
For example, any $k$-bialgebra or weak $k$-bialgebra is a very weak
bialgebra in this sense (for $\C = \Vect_k$). The structure $A$ is then
called a \emph{von Neumann core} in $\C$ if it also has an antipode
$S:A \ra A$ satisfying the axiom
\[
    \mu_3 (1 \ox S \ox 1) \delta_3 = 1:A \ra A.
\]
For example, the set of all finite paths of edges in a (row-finite)
graph algebra~\cite{IR} forms a von Neumann core in $\C = \Set$, and so does
any group in $\Set$.

Since groups $A$ in $\Set$ are characterized by the (stronger) axiom
\begin{equation}\tag{$\dagger$}\label{axiom-antipode}
    1 \ox \eta = (1 \ox \mu)(1 \ox S \ox 1)\delta_3:A \ra A \ox A,
\end{equation}
a very weak bialgebra $A$ satisfying~\eqref{axiom-antipode}, in the general
$\C$, will be called a \emph{unital} von Neumann core in $\C$. Such a unital
core $A$ always has a left inverse, namely $(1 \ox \mu)(1 \ox S \ox 1)(\delta
\ox 1)$, to the ``fusion'' operator
\[
    (1 \ox \mu)(\delta \ox 1):A \ox A \ra A \ox A,
\]
and the latter satisfies the fusion equation~\cite{S}. Any Hopf algebra in
$\C$ satisfies the axiom~\eqref{axiom-antipode}, and in this article we are
mainly interested in producing a unital von Neumann core, namely $\Endv U$,
associated to a certain type of monoidal functor $U$ into $\Vect_k$.
However, it will not be the case that all unital von Neumann cores in
$\Vect_k$ can be reproduced as such.

We will tacitly assume throughout the article that the ground
category~\cite{K} is $\Vect = \Vect_k$, for $k$ a field of characteristic
$0$, so that the categories and functors considered here are all $k$-linear
(although any reasonable category
$[\D,\Vect]$ of parameterized vector spaces would suffice). We denote by
$\Vect_f$ the full subcategory of $\Vect$ consisting of the finite
dimensional vector spaces, and we further suppose that $(\C,\ox,I,c)$ is a
braided monoidal category with a ``fibre'' functor
\[
    U:\C \ra \Vect
\]
which has both a monoidal structure $(U,r,r_0)$ and a comonoidal structure
$(U,i,i_0)$. We call $U$ \emph{separable}\footnote{Strictly, we should also
require the conditions
\begin{align*}
(r \ox 1)(1 \ox i) &= ir:UA \ox U(B \ox C) \ra U(A \ox B) \ox UC, \text{ and} \\
(1 \ox r)(i \ox 1) &= ir:U(A \ox B) \ox UC \ra UA \ox U(B \ox C)
\end{align*}
in order for $U$ to be called ``separable'', but we do not need these here.}
if $ri = 1$ and $i_0 r_0 = \dim(UI) \cdot 1$; i.e., for all $A,B \in \C$, the
diagrams
\[
\xygraph{{U(A \ox B)}="1"
    [r(2)] {UA \ox UB}="2"
    [d(1.2)] {U(A \ox B)}="3"
    "1":"2" ^-i
    "2":"3" ^-r
    "1":"3" _-1}
\qquad\qquad
\xygraph{{I}="1"
    [r(1.6)] {UI}="2"
    [d(1.2)] {I}="3"
    "1":"2" ^-{r_0}
    "2":"3" ^-{i_0}
    "1":"3" _-{\dim UI \cdot 1}}
\]
commute.

First we produce an algebra structure $(\mu,\eta)$ on
\[
    \Endv U = \int^C UC^* \ox UC
\]
using the monoidal and comonoidal structures on $U$. Secondly, we suppose
that $\C$ has a suitable small generating set $\A$ of objects, and produce a
coalgebra structure $(\delta,\epsilon)$ on $\Endv U$ when each value $UA$,
$A \in \A$, is finite dimensional. Finally, we assume that $U$ is equipped
with a natural non-degenerate form
\[
    U(A^*) \ox UA \ra k
\]
suitably related to the evaluation and coevaluation maps of $\C$ and
$\Vect_f$, where each $A \in \A$ has a $\ox$-dual $A^*$ in $\C$ which again
lies in $\A$. This last assumption is sufficient to provide $\Endv U$ with
an antipode so that it becomes a unital von Neumann core in the above sense.

By way of examples, we note that many separable monoidal functors are
constructable from separable monoidal categories; i.e., from monoidal
categories $\C$ for which the tensor product map
\[
    \ox:\C(A,B) \ox \C(C,D) \ra \C(A \ox C,B \ox D)
\]
is a naturally split epimorphism (as is the case for some finite cartesian
products such as $\Vect_f^n$). A closely related source of examples is the
notion of a weak dimension functor on $\C$ (cf.~\cite{HO}); this is a
comonoidal functor
\[
(d,i,i_0):\C \ra \Set_f
\]
for which the comonoidal transformation components
\[
i = i_{C,D}:d(C \ox D) \ra dC \times dD
\]
are injective functions, while the unique map $i_0:dI \ra 1$ is surjective.
Various examples are described at the conclusion of the paper.

We suppose the reader is familiar to some extent with the standard
references on the problem when restricted to the case of $U$ strong
monoidal.

We would like to thank Ross Street for several helpful comments.

\section{The algebraic structure on $\Endv U$}

If $\C$ is a ($k$-linear) monoidal category and
\[
    U:\C \ra \Vect
\]
has a monoidal structure $(U,r,r_0)$ and a comonoidal structure $(U,i,i_0)$,
then $\Endv U$ has an associative and unital $k$-algebra structure whose
multiplication $\mu$ is the composite map
\[
\xygraph{{\dint^C UC^* \ox UC \ox \dint^D UD^* \ox UD}="s"
    :[d(1.4)] {\dint^{C,D} UC^* \ox UD^* \ox UC \ox UD} _-\cong
    :[d(1.4)] {\dint^{C,D} (UC \ox UD)^* \ox UC \ox UD} _-{\text{can}}
    :[r(5.4)] {\dint^{C,D} U(C \ox D)^* \ox U(C \ox D)} _-{\tint i^* \ox r}
    :[u(2.8)] {\dint^B UB^* \ox UB}="t" _-{\tint^\ox} 
 "s":"t" ^-\mu}
\]
while the unit $\eta$ is given by
\[
\xygraph{{k}="s"
    :[d(1.2)] {k^* \ox k} _-\cong
    :[r(2.6)] {UI^* \ox UI.} _-{i_0^* \ox r_0}
    :[u(1.2)] {\dint^C UC^* \ox UC}="t" _-{\copr_{C = I}}
 "s":"t" ^-\eta}
\]
The associativity and unit axioms for $(\Endv U,\mu,\eta)$ now follow
directly from the corresponding associativity and unit axioms for
$(U,r,r_0)$ and $(U,i,i_0)$. An augmentation $\epsilon$ is given by
\[
\xygraph{{UD^* \ox UD}="s"
    :[u(1.4)] {\dint^C UC^* \ox UC} ^-{\copr_{C = D}}
    :[r(2)] {k}="t" ^-\epsilon
 "s":"t" _-{\mathrm{ev}}}
\]
in $\Vect$, where $\epsilon \eta = \dim UI \cdot 1$.

We also observe that the coend
\[
    \Endv U = \int^C UC^* \ox UC
\]
actually exists in $\Vect$ if $\C$ contains a small full subcategory $\A$
with the property that the family
\[
    \{Uf:UA \ra UC \mid f \in \C(A,C), A \in \A\}
\]
is epimorphic in $\Vect$ for each object $C \in \C$. In fact, we shall
use the stronger condition that the maps
\[
    \alpha_C:\int^{A \in \A} \C(A,C) \ox UA \ra UC
\]
should be isomorphisms, not just epimorphisms. This stronger condition implies
that we can effectively replace $\int^{C \in \C}$ by $\int^{A \in \A}$ since
\begin{align*}
\int^C UC^* \ox UC
    &\cong \int^C UC^* \ox (\int^A \C(A,C) \ox UA) \\
    &\cong \int^A UA^* \ox UA
\end{align*}
by the Yoneda lemma.

If we furthermore ask that each value $UA$ be finite dimensional for $A$ in
$\A$, then
\[
    \Endv U \cong \int^{A \in \A} UA^* \ox UA
\]
is canonically a $k$-coalgebra with counit the augmentation $\epsilon$, and
comultiplication $\delta$ given by
\[
\xygraph{{UA^* \ox UA}="s"
    :[u(1.5)] {\dint^A UA^* \ox UA} ^-{\copr}
    :[r(4)] {\dint^A UA^* \ox UA \ox \dint^A UA^* \ox UA}="t" ^-\delta
 "s":[r(4)] {UA^* \ox UA \ox UA^* \ox UA,} _-{1 \ox n \ox 1}
    :"t" _-{\copr\, \ox\, \copr}}
\]
where $n$ denotes coevaluation in $\Vect_f$.

\begin{proposition}
If $U$ is separable then $\Endv U$ satisfies the $k$-bialgebra axiom
\[
\xygraph{{\Endv U \ox \Endv U}="s"
    :[d(2.4)] {\Endv U} _-\mu
    :[r(3)] {\Endv U \ox \Endv U.}="t" ^-\delta
 "s":[r(3)] {(\Endv U)^{\ox 4}} ^-{\delta \ox \delta}
    :[d(1.2)] {(\Endv U)^{\ox 4}} ^-{1 \ox c \ox 1}
    :"t"  ^-{\mu \ox \mu}}
\]
\end{proposition}

\begin{proof}
Let $\B$ denote the monoidal full subcategory of $\C$ generated by $\A$ (we
will essentially replace $\C$ by this small category $\B$). Then, for all
$C,D$ in $\B$, we have, by induction on the tensor lengths of $C$ and $D$,
that $U(C \ox D)$ is finite dimensional since it is a retract of $UC \ox UD$.
Moreover, we have
\[
    \int^{A \in \A} UA^* \ox UA ~\cong~ \int^{B \in \B} UB^* \ox UB
\]
by the Yoneda lemma, since the natural family
\[
    \alpha_B:\int^{A \in \A} \C(A,B) \ox UA \ra UB
\]
is an isomorphism for all $B \in \B$. Since $ri = 1$, the triangle
\[
\xygraph{{k}="s"
    [u(0.6)r(3)] {(UC \ox UD) \ox (UC \ox UD)^*}
    :[d(1.2)] {U(C \ox D) \ox U(C \ox D)^*}="t" ^-{r \ox i^*}
 "s":[u(0.55)r(1.2)] ^-n
 "s":[d(0.55)r(1.4)] _-n}
\]
commutes in $\Vect_f$, where $n$ denotes the coevaluation maps. The asserted
bialgebra axiom then holds on $\Endv U$ since it reduces to the following
diagram on filling in the definitions of $\mu$ and $\delta$ (where, for the
moment, we have dropped the symbol ``$\ox$''):
\[
\xygraph{{UC ~ UC^* ~  UD ~ UD^*}="11"
  [d(1.2)] {UC ~ UD ~ UC^* ~ UD^*}="21"
  [d(1.2)] {UC ~ UD ~ (UC ~ UD)^*}="31"
  [d(1.2)] {U(C ~ D) ~  U(C ~ D)^*}="41"
"11"[r(5.8)] {UC ~ (UC ~ UC^*) ~ UC^* ~  UD ~ (UD ~ UD^*) ~ UD^*}="12"
  [d(1.2)] {UC ~  UD ~ UC ~ UD ~ UC^* ~ UD^* ~ UC^* ~ UD^*}="22"
  [d(1.2)] {UC ~  UD ~ UC ~ UD ~ (UC ~ UD)^* ~ (UC ~ UD)^*}="32"
  [d(1.2)] {U(C ~ D) ~ U(C ~ D) ~ U(C ~ D)^* ~ U(C ~ D)^*}="42"
"11":"21" _-\cong
"21":"31" _-\cong
"31":"41" _-{r ~ i^*}
"12":"22" ^-\cong
"22":"32" ^-\cong 
"32":"42" ^-{r ~ r ~ i^* ~ i^*}
"11":"12" ^-{1 ~ n ~ 1 ~ 1 ~ n ~ 1}
"31":"32" ^-{1 ~ n ~ 1}
"41":"42" ^-{1 ~ n ~ 1}}
\]
for all $C,D \in \B$.
\end{proof}

Notably the bialgebra axiom
\[
    \xygraph{{\Endv U \ox \Endv U}="s"
        :[r(2.4)] {\Endv U} ^-\mu
        :[l(1.2)d(1.2)] {k}="t" ^-\epsilon
    "s":"t" _-{\epsilon \ox \epsilon}}
\]
does not hold in general, while the form of the axiom
\[
    \xygraph{{k}="s"
        :[l(1.2)u(1.2)] {\Endv U} ^-\eta
        :[r(2.4)] {\Endv U \ox \Endv U}="t" ^-\delta
    "s":"t" _-{\eta \ox \eta}}
\]
holds where we multiply $\delta$ by $\dim UI$.

The $k$-bialgebra axiom established in the above proposition implies that
the ``fusion'' operator $(1 \ox \mu)(\delta \ox 1):A \ox A \ra A \ox A$
satisfies the fusion equation (see~\cite{S} for details).

The $k$-linear dual of $\Endv U$ is of course
\[
    [\int^C UC^* \ox UC,k] \cong \int_C [UC^*,UC^*]
\]
which is the endomorphism $k$-algebra of the functor
\[
    U(-)^*:\C^\op \ra \Vect.
\]
If $\ob \A$ is finite, so that
\[
    \int^A UA^* \ox UA
\]
is finite dimensional, then
\[
    \int_C [UC^*,UC^*] \cong \int_A [UA^*,UA^*]
\]
is also a $k$-coalgebra.

\section{The unital von Neumann antipode}

We now take $(\C,\ox,I,c)$ to be a braided monoidal category and $\A \subset
\C$ to be a small full subcategory of $\C$ for which the monoidal and
comonoidal functor $U:\C \ra \Vect$
induces
\[
    U:\A \ra \Vect_f
\]
on restriction to $\A$. We suppose that $\A$ is such that
\begin{itemize}
\item the identity $I$ of $\ox$ lies in $\A$, and each object of $A \in \A$
has a $\ox$-dual $A^*$ lying in $\A$.
\end{itemize}
With respect to $U$, we suppose $\A$ has the properties
\begin{itemize}
\item ``$U$-irreducibility'': $\A(A,B) \neq 0$ implies $\dim UA = \dim UB$
for all $A,B \in \A$,

\item ``$U$-density'': the canonical map
\[
    \alpha_C:\int^{A \in \A} \C(A,C) \ox UA \ra UC
\]
is an isomorphism for all $C \in \C$,

\item ``$U$-trace'': each object of $\A$ has a $U$-trace in $\C(I,I)$, where
by $U$-trace of $A \in \A$ we mean an isomorphism $d(A)$ in $\C(I,I)$ such
that the following two diagrams commute.
\[
\xygraph{{I}="s"
    :[d(1.2)] {A \ox A^*} _-n
    :[r(1.9)] {A^* \ox A} ^-c
    :[u(1.2)] {I}="t" _-e
 "s":"t" ^-{d(A)}}
\qquad\qquad
\xygraph{{k}="s"
    :[d(1.2)] {UI} _-{r_0}
    :[r(2.5)] {UI}="t" ^-{\dim UI \cdot U(d(A))}
 "s":[r(2.5)] {k} ^-{\dim UA}
    :"t" ^-{r_0}}
\]
We shall assume $\dim UI \neq 0$ so that the latter assumption implies $\dim
UA \neq 0$, for all $A \in \A$.
\end{itemize}

We require also a natural isomorphism
\[
\xygraph{{u=u_A:U(A^*)} :[r(2.1)] {UA^*} ^-\cong}
\]
such that
\begin{equation} \tag{$n,r,r_0$}
\vcenter{\xygraph{{k}="s"
    :[r(2)] {UI} ^-{r_0}
    :[d] {U(A \ox A^*)}="t" ^-{Un}
 "s":[d] {UA \ox UA^*} _-n
    :[dr] {UA \ox U(A^*)} _-{1 \ox u^{-1}}
    :"t" _-r}}
\end{equation}
commutes, and
\begin{equation} \tag{$e,i,i_0$}
\vcenter{\xygraph{{U(A^* \ox A)}="s"
    :[u] {UI} ^-{Ue}
    :[r(2)] {k}="t" ^-{i_0}
 "s":[dr] {U(A^*) \ox UA} _-i
    :[ur] {UA^* \ox UA} _-{u \ox 1}
    :"t" _-e}}
\end{equation}
commutes. This means that $U$ ``preserves duals'' when restricted to $\A$.

An endomorphism
\[
    \sigma:\Endv U \ra \Endv U
\]
may be defined by components
\[
\xygraph{{UA^* \ox UA}="s"
    :[r(3)] {U(A^*)^* \ox U(A^*),} ^-{\sigma_A}
    :[u(1.4)] {\dint^A UA^* \ox UA}="t" ^-{\copr}
 "s":[u(1.4)] {\dint^A UA^* \ox UA} _-{\copr}
    :"t" ^-\sigma}
\]
each $\sigma_A$ being given by commutativity of
\[
\xygraph{{UA^* \ox UA}="s"
    :[r(3.4)] {U(A^*)^* \ox U(A^*)}="t" ^-{\sigma_A}
 "s":[d(1.4)] {UA^* \ox UA^{**}} _-{1 \ox \rho}
    :[r(3.4)] {U(A^*) \ox U(A^*)^*} ^-{u^{-1} \ox u^*}
    :"t" _-c}
\]
where $\rho$ denotes the canonical isomorphism from a finite dimensional
vector space to its double dual. Clearly each component $\sigma_A$ is
invertible.

\begin{theorem}
Let $\C$, $\A$, and $U$ be as above, and suppose that $U$ is braided and
separable as a monoidal functor. Then there is an invertible antipode $S$ on
$\Endv U$ such that $(\Endv U, \mu,\eta,\delta,\epsilon,S)$ is a unital von
Neumann core in $\Vect_k$.
\end{theorem}

\begin{proof}
A family of maps $\{S_A \mid A \in \A\}$ is defined by
\[
    S_A = \dim UI \cdot (\dim UA)^{-1} \cdot \sigma_A.
\]
Then, by the $U$-irreducibility assumption on the category $\A$, this family
induces an invertible endomorphism $S$ on the coend
\[
    \Endv U \cong \sum_{n=1}^\infty \int^{A \in \A_n} UA^* \ox UA,
\]
where $\A_n$ is the full subcategory of $\A$ determined by $\{A \mid \dim UA
= n\}$. We now take $S$ to be the antipode on $\Endv U$ and check that
\[
    1 \ox \eta = (1 \ox \mu)(1 \ox S \ox 1)\delta_3.
\]
From the definition of $\mu$ and $\delta$, we require commutativity of the
exterior of the following diagram (where, again, we have dropped
the symbol ``$\ox$''):
\[
\xygraph{{UA^* ~  UA ~  UA^* ~  UA ~ UA^* ~ UA}="11"
    [d(3)] {UA^* ~  UA~  UA^*~  UA}="21"
    [d(2.1)] {UA^* ~ UA}="31"
    [d(2.1)] {UA^* ~ UA ~ I}="41"
"11"[r(6)] {UA^* ~ UA ~ U(A^*)^* ~ U(A^*) ~ UA^* ~ UA}="12"
    :[d(3)] {UA^* ~ UA ~ UA^{**} ~ UA^* ~ UA^* ~ UA}="22" ^-\cong
    :[d] {UA^* ~ UA ~ U(A^*)^* ~ U(A^*) ~ UA^* ~ UA} ^-\cong
    :[d] {UA^* ~ UA ~ (U(A^*) ~ UA)^* ~ U(A^*) ~ UA} ^-\cong
    :[d] {UA^* ~ UA ~ U(A^* ~ A)^* ~  U(A^* ~ A)} ^-{1~ 1~ i^*~ r}
    :[d(1.2)] {UA^* ~ UA ~ \dint^B UB^* ~ UB}="42" ^-{1~ 1~ \copr}
"21":[r(2.5)u] {UA^* ~ UA ~ UA ~ UA^* ~ UA^* ~ UA} _-{1 ~ 1 ~ n ~ 1 ~ 1}
    :[u] {UA^* ~ UA ~ UA ~ UA^* ~ UA^* ~ UA} _-{1 ~ 1 ~ 1 ~ c ~ 1}
    :"11" _-{1 ~ 1 ~ c ~ 1 ~ 1}
    "21":"11" ^-{1 ~ n ~ 1 ~ 1 ~ 1}
    "31":"21" ^-{1 ~ n ~ 1}
    "41":"31" ^-\cong
    "11":"12" ^-{1 ~ 1 ~ S_A ~ 1 ~ 1}
    "21":"22" ^-{1 ~ 1 ~ e^* ~ 1 ~ 1}
    "41":"42" ^-{1 ~ 1 ~ \eta}
    "11"[d(1.5)r(0.7)] {\scriptstyle (1)}
    "11"[d(1.5)r(5)] {\scriptstyle (2)}
    "21"[d(2.1)r(2.5)] {\scriptstyle (3)}}
\]
The region labelled by $(1)$ commutes on composition with $1 \ox n \ox 1$
since
\[
\xygraph{{k}="11"
    :[r(3.5)] {UA \ox  UA^*}="12" ^-n
    :[d(1.2)] {UA \ox UA \ox UA^* \ox UA^*}="13" ^-{1 \ox n \ox 1}
    :[d(1.2)] {UA \ox UA \ox UA^* \ox UA^*}="14" ^-{1 \ox 1 \ox c}
    :[d(1.2)] {UA \ox UA^* \ox UA \ox UA^*}="15" ^-{1 \ox c \ox 1}
"11":[d(3.6)] {UA \ox UA^*} _-n
    :"15" ^-{n \ox 1 \ox 1}}
\]
commutes (choose a basis for $UA$). The region labelled by $(2)$ now
commutes by inspection of:
\[
\xygraph{{UA \ox UA^*}="11"
    :[u(1.33)] {UA \ox UA \ox UA^* \ox UA^*}="21" ^-{1 \ox n \ox 1}
    :[u(1.33)] {UA \ox UA \ox UA^* \ox UA^*}="31" ^-{1 \ox 1 \ox c}
    :[u(1.33)] {UA \ox UA^* \ox UA \ox UA^*}="41" ^-{1 \ox c \ox 1}
"11":[r(5)] {UA \ox UA^{**} \ox UA^* \ox UA^*}="12" ^-{1 \ox e^* \ox 1}
    :[u(0.99)] {UA \ox UA^{**} \ox UA^* \ox UA^*}="22" _-{1 \ox 1 \ox c}
    :[u] {UA \ox UA^* \ox UA^{**} \ox UA^*}="32" _-{1 \ox c \ox 1}
    :[u] {UA \ox U(A^*) \ox U(A^*)^* \ox UA^*}="42"
    :[u] {UA \ox U(A^*)^* \ox U(A^*) \ox UA^*}="52"
"32":"42" _-{1 \ox u^{-1} \ox u^* \ox 1}
"42":"52" _-{1 \ox c \ox 1}
"41":"52" ^-{1 \ox \sigma_A \ox 1}
"41":"32" |-{1 \ox 1 \ox \rho \ox 1}
"21":"12" |-{1 \ox \rho \ox 1 \ox 1}}
\]
where the top leg of $(2)$ has been rescaled by a factor of $(\dim UI)^{-1}
\cdot \dim UA$.

From the definition of the $U$-trace $d(A)$ of $A \in \A$, we have that
\[
\xygraph{{k}="s"
    :[r(2.8)] {k} ^-{\dim UI \cdot (\dim UA)^{-1}}
    :[d(1.2)] {UI}="t" ^-{r_0}
 "s":[d(1.2)] {UI} _-{r_0}
    :"t" _-{U(d(A)^{-1})}}
\]
commutes, so that the exterior of
\[
\xygraph{{k}="a"
    :[u(1.2)] {k}="b" ^-{\dim UI \cdot (\dim UA)^{-1}}
    :[u(1.2)] {UA \ox UA^*} ^-n
    :[r(1.5)u] {UA \ox UA^*}="t" ^-{1 \ox u^{-1}}
    :[r(1.5)d] {U(A \ox A^*)}="c" ^-r
 "a":[r(3)] {UI} ^-{r_0}
    :[u(1.2)] {UI}="d" _-{U(d(A)^{-1})}
    :"c" _-{Un}
 "b":"d" _-{r_0}
 "t"[d(1.2)] {\scriptstyle (n,r,r_0)}}
\]
commutes.

Thus the region labelled by $(3)$, with the top leg rescaled by the
factor $\dim UI \cdot (\dim UA)^{-1}$, commutes on examination of the
following diagram:
\[
\xygraph{{k^* \ox k}="0"
    [u(3.4)] {k^* \ox UA \ox UA^*}="1"
    [u(2.8)] {k^* \ox UA^* \ox UA}="2"
    [r(3.2)u] {(UA^* \ox UA)^* \ox UA^* \ox UA}="3"
    [r(3.2)d] {(U(A^*) \ox UA)^* \ox U(A^*) \ox UA}="4"
    [d(3)] {U(A^* \ox A)^* \ox U(A^* \ox A)}="5"
    [d(3.2)] {\dint^B UB^* \ox UB}="6"
 "2":[d(1)r(1.8)] {k^* \ox U(A^*) \ox UA}="a1" |-{1 \ox u^{-1} \ox 1}
    :[d(1)r(1.8)] {UI^* \ox U(A^* \ox A)}="a2" |-{i^*_0 \ox r}
    :"5" ^-{Ue^* \ox 1}
 "1":[d(1)r(1.8)] {k^* \ox UA \ox U(A^*)}="b1" |-{1 \ox 1 \ox u^{-1}}
    :[d(1)r(1.8)] {UI^* \ox U(A \ox A^*)}="b2" |-{i^*_0 \ox r}
 "0"[r(3.6)] {UI^* \ox UI}="z"
    [r(2)u(1.2)] {UI^* \ox UI}="w"
"z":"w" ^-1
"w":"6" _-\copr
"0":"z" _-{i^*_0 \ox r}
"z":"6" _-{\copr}
"0":"1" |-{1 \ox \dim UI \cdot (\dim UA)^{-1} \cdot n}
"1":"2" ^-{1 \ox c}
"2":"3" ^-{e^* \ox 1 \ox 1\quad}
"3":"4" ^-{\qquad(u \ox 1)^* \ox (u^{-1} \ox 1)}
"4":"5" ^-{i^* \ox r}
"5":"6" ^-{\copr}
"b1":"a1" ^-{1 \ox c}
"b2":"a2" ^-{1 \ox Uc}
"a2":"w" ^<(0.6){1 \ox Ue}
 "z":"b2" ^-{1 \ox U(d(A)^{-1}) \cdot Un}
 "0"[ur] {\scriptstyle (n,r,r_0)}
 "b2"[u(1.8)l(0.9)] {\scriptstyle (*)}
 "3"[d] {\scriptstyle (e,i,i_0)}}
\]
whose commutativity depends on the hypothesis that $(U,r,r_0)$ is braided
monoidal in order for
\[
\xygraph{{UA \ox U(A^*)}="s"
    :[r(2.6)] {U(A^*) \ox UA} ^-c
    :[d(1.2)] {U(A^* \ox A)}="t" ^-r
 "s":[d(1.2)] {U(A \ox A^*)} _-r
    :"t" _-{Uc}
 "s"[r(1.2)d(0.6)]{\scriptstyle (*)}}
\]
to commute.
\end{proof}

\section{The fusion operator}

Let $E = \Endv U$. The unital von Neumann axiom on $E$ implies that the fusion
operator
\[
    f = (1 \ox \mu)(\delta \ox 1):E \ox E \ra E \ox E
\]
has a left inverse, namely $g=(1 \ox \mu)(1 \ox S \ox 1)(\delta \ox 1)$. For
this we consider the following diagram:
\[
\xygraph{{E \ox E}="11"
    [r(2)] {E^{\ox 3}}="12"
    [r(2)] {E \ox E}="13"
"12"[d(1.2)] {E^{\ox 4}}="22"
    [r(2)d(0.6)] {E^{\ox 3}}="33"
    [l(2)d(0.6)] {E^{\ox 4}}="42"
    [l(2)] {E^{\ox 3}}="41"
    [d(1.2)] {E \ox E}="51"
    [r(4)] {E^{\ox 3}.}="53"
"11":"12" ^-{\delta \ox 1}
"11":"41" ^<(0.57){1 \ox \eta \ox 1}
"11":"22" _-{\delta_3 \ox 1}
"11":@/_5ex/"51" _-1
"12":"13" ^-{1 \ox \mu}
"12":@<-3pt>"22" _-{1 \ox \delta \ox 1}
"12":@<3pt>"22" ^-{\delta \ox 1 \ox 1}
"13":"33" ^-{\delta \ox 1}
"22":"33" ^-{1 \ox 1 \ox \mu}
"22":"42" ^-{1 \ox S \ox 1 \ox 1}
"33":"53" ^-{1 \ox S \ox 1}
"41":"51" ^-{1 \ox \mu}
"42":"41" _-{1 \ox \mu \ox 1}
"42":"53" ^-{1 \ox 1 \ox \mu}
"53":"51" _-{1 \ox \mu}}
\]
In particular $f=(1 \ox \mu)(\delta \ox 1)$ is a partial isomorphism, i.e.,
$fgf=f$ and $gfg=g$.

\section{Examples of separable monoidal functors in the present context}

Unless otherwise indicated, categories, functors, and natural
transformations shall be $k$-linear, for $k$ a suitable field. 

For these examples we recall that a (small) $k$-linear promonoidal category
$(\A,p,j)$ (previously called ``premonoidal'' in~\cite{PC}) consists of a
$k$-linear category $\A$ and two $k$-linear functors
\begin{align*}
p &: \A^\op \ox \A^\op \ox \A \ra \Vect \\
j &: \A \ra \Vect
\end{align*}
equipped with associativity and unit constraints satisfying axioms (as
described in~\cite{PC}) analogous to those used to define a monoidal
structure on $\A$. The notion of a symmetric promonoidal category
(also introduced in~\cite{PC}) was extended in~\cite{DPS} to that of a
braided promonoidal category.

The main point is that (braided) promonoidal structures on $\A$ correspond
to cocontinuous (braided) monoidal structures on the functor category
$[\A,\Vect]$. This latter monoidal structure is often called the convolution
product of $\A$ and $\Vect$.

\begin{example}
Let $(\A,p,j)$ be a small braided promonoidal category with 
\[
    \A(I,I) \cong I = k \qquad \text{and} \qquad j=\A(I,-),
\]
and suppose that each hom-space $\A(a,b)$ is finite
dimensional. Let $f \in [\A,\Vect_f]$ be a very weak bialgebra in the
convolution $[\A,\Vect]$. Suppose also that $\A \subset \C$ where $\C$ is a
separable braided monoidal category with
\[
    p(a,b,c) \cong \C(a \ox b,c)
\]
naturally; we suppose the induced maps
\begin{equation}\tag{$\ddagger$}\label{pc}
    \int^c p(a,b,c) \ox \C(c,C) \ra \C(a \ox b,C)
\end{equation}
are isomorphisms (e.g., $\A$ monoidal). We also suppose that each $a \in \A$
has a dual $a^* \in \A$. Then we have maps
\[
    \mu : f * f \ra f \qquad \text{and} \qquad \eta : k \ra fI
\]
and
\[
    \delta : f \ra f * f \qquad \text{and} \qquad \epsilon : fI \ra k
\]
satisfying associativity and unital axioms.

Define the functor $U:\C \ra \Vect$ by
\[
    U(C) = \int^a fa \ox \C(a,C);
\]
then, by the Yoneda lemma, $U(a^*) \cong U(a)^*$ if $f(a^*) \cong f(a)^*$
for $a \in \A$. Moreover, $U$ is monoidal and comonoidal on $\C$ via the maps
$r$ and $i$ described in the diagram:
\[
\xygraph{{UC \ox UD}="x"
    :[r(3.8)] {\dint^{a,b} fa \ox fb \ox \C(a,C) \ox \C(b,D)}="1" ^-\cong
    [d(1.6)] {\dint^{a,b} fa \ox fb \ox \C(a \ox b,C \ox D)}="2"
    [d(1.6)] {\dint^{a,b} fa \ox fb \ox \dint^c p(a,b,c) \ox \C(c,C \ox D)}="3"
    [d(1.6)] {\dint^c fc \ox \C(c,C \ox D),}="4"
    :[l(3.8)] {U(C \ox D)}="y" _-=
"x":@<-4pt>"y" _-r
"y":@<-4pt>"x" _-i
"1":@<-3pt>"2"
"2":@<-3pt>"1" _-{\textrm{\quad $\C$ separable}}
"2":"3" ^-{\quad \eqref{pc}}
"3":@<-3pt>"4" _-\mu
"4":@<-3pt>"3" _-\delta}
\]
Thus, if $f$ is separable, then so is $U$ with $\dim UI = \dim fI$ since
\[
    UI = \int^a fa \ox \C(a,I) \cong fI
\]
by the Yoneda lemma, so that $i_0 r_0 = \dim UI \cdot 1$ if and only if
$\epsilon \eta = \dim fI \cdot 1$.
\end{example}

\begin{example}
Suppose that $(\A^\op,p,j)$ is a small promonoidal category with $I \in \A$
such that $j \cong \A(-,I)$ and with each $x \in \A$ an ``atom'' in $\C$
(i.e., an object $x \in \C$ for which $\C(x,-)$ preserves all colimits)
where $\C$ is a cocomplete and cocontinuous braided monoidal category
containing $\A$ and each $x \in \A$ has a dual $x^* \in \A$. Suppose that
the inclusion $\A \subset \C$ is dense over $\Vect$ (that is, the canonical
evaluation morphism
\[
    \int^a \C(a,C) \cdot a \ra C
\]
is an isomorphism for all $C \in \C$) and
\[
    x \ox y \cong \int^z p(x,y,z) \cdot z
\]
so that
\begin{align*}
\C(a,x \ox y) &= \C(a,\int^z p(x,y,z) \cdot z) \\
              &\cong \int^z p(x,y,z) \ox \C(a,z) \quad \text{since $a \in \A$
                is an atom in $\C$,} \\
              &\cong p(x,y,a) \qquad \text{by the Yoneda lemma applied to $z
                \in \A$.}
\end{align*}

Let $W:\A \ra \Vect$ be a strong promonoidal functor on $\A$. This means
that we have structure isomorphisms
\begin{align*}
    Wx \ox Wy &\cong \int^z \C(z,x \ox y) \ox Wz \\
    k &\cong WI
\end{align*}
satisfying suitable associativity and unital coherence axioms. Define the
functor $U:\C \ra \Vect$ by
\[
    UC = \int^a \C(a,C) \ox Wa.
\]
Then
\begin{align*}
    U(x^*) &= \int^a \C(a,x^*) \ox Wa \\
           &\cong W(x^*) \\
           &\cong W(x)^*,
\end{align*}
if $W(x^*) \cong W(x)^*$ for all $x \in \A$, and
\begin{align*}
    UI &= \int^a \C(a,I) \ox Wa \\
           &\cong WI \\
           &\cong k,
\end{align*}
so that $i_0 r_0 = 1$ and $r_0 i_0 = 1$. Also there are mutually inverse
composite maps $r$ and $i$ given by:
\begin{align*}
r:UC \ox UD &\cong \int^{x,y} \C(x,C) \ox \C(y,D) \ox Ux \ox Uy \\
            &\cong \int^{x,y} \C(x,C) \ox \C(y,D) \ox Wx \ox Wy \\
            &\cong \int^{x,y} \C(x,C)\ox \C(y,D)\ox \int^z \C(z,x \ox y) \ox Wz \\
            &\cong \int^z \C(z,C \ox D) \ox Wz \\
            &\cong U(C \ox D),
\end{align*}
which uses the assumptions that $\C$ is cocontinuous monoidal and $\A
\subset \C$ is dense. Thus $ri = 1$ and $ir = 1$ so that $U$ is a strong
monoidal functor.
\end{example}

\begin{example}
(See~\cite{HO} Proposition~3.) Let $\C$ be a braided compact monoidal category
and let $\A \subset \C$ be a full finite discrete Cauchy generator of $\C$
which contains $I$ and is closed under dualization in $\C$. As in the
H\"aring-Oldenburg case~\cite{HO}, we suppose that each hom-space $\C(C,D)$
is finite dimensional with a chosen natural isomorphism $\C(C^*,D^*) \cong
\C(C,D)^*$.

Then we have a separable monoidal functor
\[
    UC = \bigoplus_{a,b \in \A} \C(a,C \ox b),
\]
whose structure maps are given by the composites
\begin{align*}
r:UC \ox UD \cong&~ \bigoplus_{a,b,c,d} \C(c,C \ox b) \ox \C(a,D \ox d) \\
    \xymatrix@C=30pt@1{\ar@<2pt>[r]^-{\textrm{adjoint}} & \ar@<2pt>[l]^-{c=d}}&~
                \bigoplus_{a,b,c} \C(c,C \ox b) \ox \C(a,D \ox c) \\
        \cong&~ \bigoplus_{a,b} \C(a,D \ox (C \ox b)) \\
        \cong&~ \bigoplus_{a,b} \C(a,(D \ox C) \ox b) \\
        \cong&~ \bigoplus_{a,b} \C(a,(C \ox D) \ox b) \\
            =&~ U(C \ox D),
\end{align*}
and $r_0:k \ra UI$ the diagonal, with $i_0$ its adjoint. Also
\begin{align*}
U(C^*) &= \bigoplus_{a,b} \C(a,C^* \ox b) \\
            &\cong \bigoplus_{a,b} \C(a^*,C^* \ox b^*) \\
            &\cong \bigoplus_{a,b} \C(a,C \ox b)^* \\
            &\cong UC^*
\end{align*}
for all $C \in \C$.
\end{example}

\begin{example}
Let $(\A,p,j)$ be a finite braided promonoidal category over $\Set_f$ with
$I \in \A$ such that $j \cong \A(I,-)$ and with a promonoidal functor
\[
    d:\A^\op \ra \Set_f
\]
for which each structure map
\[
    u:\int^z p(x,y,z) \times dz \ra dx \times dy
\]
is an injection, and $u_0:dI \ra 1$ is a surjection. Then we have
corresponding maps
\[
    \dint^z k[p(x,y,z)] \ox k[dz]~~
        \xy
            \ar@{>->} (0,1); (7,1)
            \ar@{->>} (7,-1); (-1.5,-1)
        \endxy
    ~~k[dx] \ox k[dy]
\]
and
\[
    k[dI]~~
        \xy
            \ar@{->>} (0,1); (8.5,1)
            \ar@{>->} (7,-1); (0,-1)
        \endxy
    ~~k[1],
\]
where $k[s]$ denotes the free $k$-vector space on the (finite) set $s$, in
$\Vect_f$. Define the functor $U:\C \ra \Vect_f$ by
\[
    Uf = \int^x fx \ox k[dx]
\]
for $f \in \C = [k_* \A,\Vect_f]$ (with the convolution braided monoidal
closed structure) so that
\begin{align*}
r:Uf \ox Ug =&~ \left(\int^x fx \ox k[dx]\right) \ox \left(\int^y gx \ox
               k[dy]\right) \\
    \cong&~ \int^{x,y} fx \ox gy \ox (k[dx] \ox k[dy]) \\
    \xymatrix@1{\ar@<2pt>[r] & \ar@<2pt>[l]}&~
                \int^{x,y} fx \ox gy \ox \left(\int^z k[p(x,y,z)] \ox k[dz]\right) \\
    \cong&~ \int^z \left( \int^{x,y} fx \ox gy \ox k[p(x,y,z)]\right) \ox k[dz] \\
    =&~ \int^z (f \ox g)(z) \ox k[dz] \\
    =&~ \int^z U(f \ox g)
\end{align*}
and
\begin{align*}
    i_0:UI =&~ \int^x k[\A(I,x)] \ox k[dx] \\
        \cong&~ k[dI] \\
        \xymatrix@1{\ar@<2pt>[r] & \ar@<2pt>[l]}&~ k[1] \cong k.
\end{align*}
Hence $i_0 r_0 = \dim UI \cdot 1 = |dI| \cdot 1$. Thus, $U$ becomes a
separable monoidal functor.
\end{example}

\begin{example}
Let $\A$ be a finite (discrete) set and give the cartesian product $\A
\times \A$ the $\Set_f$-promonoidal structure corresponding to bimodule
composition (i.e., to matrix multiplication). If
\[
    d:\A \times \A \ra \Set_f
\]
is a promonoidal functor, then its associated structure maps
\begin{align*}
\sum_{z,z'} p((x,x'),(y,y'),(z,z')) \times d(z,z')
    =&~ \sum_{z,z'} \A(z,x) \times \A(x',y) \times \A(y',z') \times d(z,z') \\
\cong&~ \A(x',y) \times d(x,y') \\
 \ra&~ d(x,x') \times d(y,y'),
\end{align*}
and
\begin{align*}
    \sum_{z,z'} j(z,z') \times d(z,z')
        =&~ \sum_{z,z'} \A(z,z') \times d(z,z') \\
    \cong&~ \sum_z d(z,z) \\
    \ra&~ 1,
\end{align*}
are determined by components
\begin{align*}
    d(x,y') &~\xymatrix@C=20pt{\ar@{>->}[r]&} d(x,y) \times d(y,y')  \\
    d(z,z)  &\xymatrix{\ar@{->>}[r]&} 1
\end{align*}
which give $\A$ the structure of a discrete cocategory over $\Set_f$.

Define the functor $U:\C = [k_* (\A \times \A),\Vect_f] \ra \Vect_f$ by
\[
    Uf = \bigoplus_{x,y} (f(x,y) \ox k[d(x,y)]).
\]
Then we obtain monoidal and comonoidal structure maps
\begin{align*}
U(f \ox g) &\xymatrix@C=30pt{\ar@<-2pt>[r]_i & \ar@<-2pt>[l]_r} Uf \ox Ug \\
UI &\xymatrix@C=30pt{\ar@<-2pt>[r]_{i_0} & \ar@<-2pt>[l]_{r_0}} k \cong k[1]
\end{align*}
from the canonical maps
\begin{multline*}
    \qquad
    \bigoplus_{x,y,z} f(x,z) \ox g(z,y) \ox k[d(x,y)] \\
    \xymatrix@C=40pt{\ar@<-2pt>[r]_{z=u=v} & \ar@<-2pt>[l]_{\textrm{adjoint}}}
    \bigoplus_{x,u} \big(f(x,u) \ox k[d(x,u)]\big) \ox 
    \bigoplus_{v,y} \big(g(v,y) \ox k[d(v,y)]\big)
    \qquad
\end{multline*}
and
\[
    \bigoplus_z k[d(z,z)]
    \xymatrix@C=20pt{\ar@<-2pt>[r] & \ar@<-2pt>[l]}
    k \cong k[1].
\]
These give $U$ the structure of a separable monoidal functor on $\C$.
\end{example}

\section{Concluding remarks}

If the original ``fibre'' functor $U$ is faithful and exact then the Tannaka
equivalence (duality)
\[
    \Lex(\C^\op,\Vect) \simeq \Comod(\Endv U)
\]
is available. Thus, since $\C$ is braided monoidal, so is $\Comod(\Endv U)$
with the tensor product and unit induced by the convolution product on
$\Lex(\C^\op,\Vect)$; for convenience we recall~\cite{D} that, for $\C$
compact, this convolution product is given by the restriction to
$\Lex(\C^\op,\Vect)$ of the coend
\begin{align*}
    F * G &= \int^{C,D} FC \ox GD \ox \C(-,C \ox D) \\
          &\cong \int^C FC \ox G(C^* \ox -)
\end{align*}
computed in the whole functor category $[\C^\op,\Vect]$. Moreover, when $U$
is separable monoidal, the category $\Co(\Endv U)$ of cofree coactions of
$\Endv U$ (as constructed in~\cite{JS} for example) also has a monoidal
structure $(\Co(\Endv U),\ox,k)$, this time obtained from the algebra
structure of $\Endv U$. The forgetful inclusion
\[
   \Comod(\Endv U) \subset \Co(\Endv U)
\] 
preserves colimits while $\Comod(\Endv U)$ has a small generator, namely
$\{UC \mid C \in \C\}$, and thus, from the special adjoint functor theorem,
this inclusion has a right adjoint. The value of the adjunction's counit at
the functor $F \ox G$ in $\Co(\Endv U)$ is then a split monomorphism and, in
particular, the monoidal forgetful functor
\[
    \Comod(\Endv U) \ra \Vect,
\]
which is the composite $\Comod(\Endv U) \subset \Co(\Endv U) \ra \Vect$,
is a separable monoidal functor extension of the given functor $U:\C \ra
\Vect$.


\end{document}